\documentclass[12pt]{amsart}

\begin{document}
\title[Fonctions holomorphes injectives.] {Une caract\'erisation des fonctions holomorphes injectives en analyse ultram\'etrique.}
%%%\subjclass{30C65}
%%%\keywords{}

\author[J. Rivera-Letelier]{Juan Rivera-Letelier}
\address{J. Rivera-Letelier \\
         Mathematics Department \\
         SUNY at Stony Brook \\
         Stony Brook, NY 11794-3660}
\email{rivera@math.sunysb.edu}
%\thanks{}

\date{Avril 2002}

%\begin{abstract}
%\end{abstract}

\maketitle

%\setcounter{tocdepth}{1}
%\tableofcontents

%\setcounter{page}{1}
%\renewcommand{\thepage}{\arabic{page}}

%\usepackage{amsmath,amssymb}
%\usepackage{amsfonts}
%\usepackage{graphicx}
%\usepackage{epsfig}
%\usepackage{fullpage}
%\usepackage{makeidx}
%\usepackage{latin1}{inputenc}
%\usepackage{french}
%\font\fiverm=cmr10 scaled 600  
%\input amssym.def

%\bibliographystyle{amsalpha}

%\input prepictex.tex  
%\input pictex.tex  
%\pagestyle{plain}  
%\setlength{\textwidth}{13.5cm}  
%\setlength{\textheight}{22cm}  
%\oddsidemargin 20pt  
%\topmargin -1cm  
%\parindent 25pt  
%\parskip 1ex
%\addtolength\topmargin{-\headheight}  
%\addtolength\topmargin{-\headsep} 

%\theoremstyle{plain} 
\newtheorem{proposition}{Proposition}
\newtheorem{lemme}{Lemme}
\newtheorem{remarque}{Remarque}
\newtheorem{definition}{D\'efinition}  
\newtheorem{corollaire}{Corollaire}  
\newtheorem{exemple}{Exemple}
\newtheorem{theoreme}{Th\'eor\`eme}

\newcommand{\diam}{\mbox{\rm diam}}
\newcommand{\dist}{\mbox{\rm dist}}
\newcommand{\can}{{\mbox{\rm can}}}
\renewcommand{\mod}{\, \, \mbox{\rm mod}}

\def\tl{\tilde}
\def\ov{\overline}
\def\cal{\mathcal}

\def\e{\varepsilon}
\def\la{\lambda}

\def\CC{\Bbb C}
\def\DD{\Bbb D}
\def\FF{\Bbb F}
\def\HH{\Bbb H}
\def\PP{\Bbb P}
\def\QQ{\Bbb Q}
\def\RR{\Bbb R}
\def\TT{\Bbb T}
\def\ZZ{\Bbb Z}

\def\preuve{\par\noindent {\bf Preuve.} }
\def\pp{{\PP(\CC_p)}}
\def\f{{\ov{\F}_p}}
\def\pf{{\PP(\ov{\FF}_p)}}

\def\cf{{\cal F}}	
\def\ci{{\cal I}}
\def\co{{\cal O}}
\def\cp{{\cal P}}
\def\cs{{\cal S}}
\def\sc{{\cal S}}
\def\ct{{\cal T}}

\def\ha{\widehat{\cal A}}
\def\he{\widehat{\cal E}}
\def\hr{\widehat{\cal R}}	
\def\hv{\widehat{\cal V}}	
\def\hw{\widehat{\cal W}}

\begin{quote}
{\footnotesize
\noindent
{\bf R\'esum\'e.}
	On montre qu'une fonction holomorphe non-constante $f$ d\'efinie sur un sous-espace analytique de $\CC_p$ est injective si et seulement si on a
$$
	\left| \frac{f(x) - f(y)}{x - y} \right|^2 = 
		|f'(x) \cdot f'(y)|, 
	\mbox{ pour tous $x$ et $y$ distincts.}
$$
	Cette caract\'erisation d\'emontre l'analogue, pour les fonctions holomorphes, d'une conjecture de A.~Escassut et M.C.~Sarmant.
	D'autre part on donne une contre-exemple \`a cette conjecture, qui concerne les \'el\'ements bi-analytiques.
}

\

{\footnotesize
\noindent
{\bf Abstract.}
	We prove that a non constant holomorphic function $f$ defined over an analytic subspace of $\CC_p$ is injective if and only if
$$
	\left| \frac{f(x) - f(y)}{x - y} \right|^2 = 
		|f'(x) \cdot f'(y)|, 
	\mbox{ for every distinct $x$ and $y$.}
$$
	This caracterisation proves the analogue, for holomorphic functions, of a conjecture of A.~Escassut and M.C.~Sarmant.
	On the other hand we give a counter example to this conjecture, that concerns bi-analytic elements.
}
\end{quote}

\

\

	Fixons un nombre premier $p$ et soient $\QQ_p$ le corps des nombres $p$-adiques et $\CC_p$ la plus petite extension alg\'ebriquement close  et compl\`ete de $\QQ_p$.
	On note $\CC_p^*$ le groupe multiplicatif des \'el\'ements non-nuls de $\CC_p$ et on note $|\CC_p^*| = \{ |z| \mid z \in \CC_p^* \}$, o\`u $| \cdot |$ est la norme de $\CC_p$.

	Une {\it boule ferm\'ee} (resp. {\it ouverte}) de $\CC_p$ est un ensemble de la forme
$$
	\{ z \in \CC_p \mid |z - a| \le r \}
		\mbox{ (resp. $\{ z \in \CC_p \mid |z - a| < r \}$)},
$$
o\`u $a \in \CC_p$ et $r \in |\CC_p^*|$.
	Un {\it affino\"{\i}de} ({\it ferm\'e}) est un ensemble de la forme $B - (B_1 \cup \ldots \cup B_n)$, o\`u $B$ est une boule ferm\'ee et les boules $B_1, \ldots, B_n$ sont ouvertes.
	Un {\it espace analytique} est une r\'eunion croissante d'affino\"{\i}des.

	Une {\it fonction holomorphe d\'efinie sur un affino\"{\i}de} est une limite uniforme de fonctions rationnelles sans p\^oles dans l'affino\"{\i}de.
	Une {\it fonction holomorphe} est une fonction d\'efinie sur un espace analytique $X \subset \CC_p$, telle que sa restriction \`a tout affino\"{\i}de contenu dans $X$ est holomorphe, voir e.g. \cite{FvP} ou \cite{Y}.

\begin{theoreme}\label{injectivite}
	Soit $X \subset \CC_p$ un espace analytique et soit $f$ une fonction holomorphe non-constante.
	Alors $f$ est injective si et seulement si on a
\begin{equation}\label{eqn}
	\left| \frac{f(x) - f(y)}{x - y} \right|^2 = 
		|f'(x) \cdot f'(y)|, 
	\mbox{ pour tous } x, y \in X \mbox{ distincts.}
\end{equation}
\end{theoreme}
	
	Pour $a, b, c, d \in \CC_p$ distincts on appelle
$
	R(a, b~; c, d) = \frac{(a-c)(b-d)}{(a-d)(b-c)} \in \CC_p^*
$
le {\it birapport} de $(a, b, c, d)$.
	La preuve du corollaire suivant est ci-dessous.
\begin{corollaire}\label{corollaire}
	Soit $f$ une fonction holomorphe et injective d\'efinie sur un espace analytique $X \subset \CC_p$.
	Alors pour tous $a, b, c, d \in \CC_p$ distincts, le birapport $R(f(a), f(b)~; f(c), f(d))$ est bien d\'efini et on a
$$
	|R(a, b~; c, d)| = |R(f(a), f(b)~; f(c), f(d))|.
$$
\end{corollaire}

	Il est facile de voir que toute fonction qui pr\'eserve le birapport co\"{\i}ncide avec une homographie.
	Donc on peut interpr\'eter le Th\'eor\`eme~\ref{injectivite} en disant que les fonctions holomorphes qui sont injectives sont celles qui sont ``proches'' des homographies.

\

	La d\'emonstration que la propri\'et\'e ($\ref{eqn}$) implique que $f$ est injective est simple et assez g\'en\'erale.
	Pour montrer l'implication inverse, on consid\`ere l'application $f_* : \widehat{X} \longrightarrow \HH_p$ induite par $f$, o\`u $\HH_p$ est l'espace hyperbolique $p$-adique et $\widehat{X}$ est l'enveloppe convexe de $X$ dans $\HH_p$~; voir \cite{these}, \cite{Y}, \cite{hyp}. 
	La propri\'et\'e ($\ref{eqn}$) suit facilement du fait que, quand $f$ est injective, l'application $f_*$ induit une isom\'etrie entre $\widehat{X}$ et $f_*(\widehat{X})$, voir \cite{these}, \cite{Y}.
\begin{proof}[Preuve du Corollaire~\ref{corollaire}.]
	Comme $f$ est injective, les points $f(a)$, $f(b)$, $f(c)$, $f(d) \in \CC_p$ sont distincts et donc le birapport $R(f(a), f(b)~; f(c), f(d))$ est bien d\'efini.
	Si l'on applique la propri\'et\'e ($\ref{eqn}$) aux paires $(a, c)$, $(b, d)$, $(a, d)$ et $(b, c)$ on obtient
$$
	\frac{|R(f(a), f(b)~; f(c), f(d))|^2}{|R(a, b~; c, d)|^2}
		= \frac{|f'(a) \cdot f'(c)| \cdot |f'(b) \cdot f'(d)|}
			{|f'(a) \cdot f'(d)| \cdot |f'(b) \cdot f'(c)|} = 1.
$$
\end{proof}
%%%%%%%%%%%%%%%%%%%%%%%%%%%%%%%%%%%%%%%%%%%%
\section{Sur les \'el\'ements analytiques et bi-analytiques.}
	Pour des r\'ef\'erences g\'en\'erales \`a cette section on renvoie le lecteur \`a \cite{Es95}.

	Rappelons qu'un {\it \'el\'ement analytique} d\'efini sur une partie $D \subset \CC_p$, est une limite uniforme de fonctions rationnelles sans p\^oles dans $D$.
	Un \'el\'ement analytique injectif dont l'inverse est aussi un \'el\'ement analytique, est appel\'e {\it bi-analytique}.
	Rappelons finalement qu'une partie $D \subset \CC_p$ est {\it infraconnexe} si pour tous $a \in \CC_p$ et $r' > r> 0$ tels que la couronne $\{ r < |z - a| < r' \}$ soit disjointe de $D$, on a $D \subset \{ |z - a| \le r \}$ ou $D \cap \{ |z- a| < r' \} = \emptyset$.

\

	L'\'equivalence du Th\'eor\`eme~\ref{injectivite} n'est pas valable dans d'autres contextes.
	En effet A.~Escassut et M.C.~Sarmant ont d\'ecrit un ouvert infraconnexe $D \subset \CC_p$ o\`u le polyn\^ome $P_0(z) = z^2$ est injectif et tel que la fonction $P_0|_D$ ne satisfait pas la propri\'et\'e ($\ref{eqn}$), voir \cite{ES} p.~161.
	De plus, ils ont remarqu\'e que l'\'el\'ement analytique $P_0|_D$ n'est pas bi-analytique.

	Motiv\'es par cet exemple, A.~Escassut et M.C. Sarmant ont conjectur\'e que {\it pour un \'el\'ement analytique d\'efini sur un ouvert infraconnexe, la propri\'et\'e} ($\ref{eqn}$) {\it est \'equivalente \`a ce que l'\'el\'ement analytique soit bi-analytique}, voir Conjecture~3 de \cite{ES} (p.~162).

	Dans \cite{elements} on montre une partie de cette conjecture~: tout \'el\'ement bi-analytique d\'efini sur un ouvert infraconnexe satisfait la propri\'et\'e ($\ref{eqn}$).
	D'autre part l'exemple~\ref{exemple} ci-dessous est un contre-exemple \`a l'implication inverse.

\

	On remarque que l'inverse d'une fonction holomorphe et injective est aussi une fonction holomorphe, voir e.g. \cite{Y} ou \cite{Mo}.
	Donc le Th\'eor\`eme~\ref{injectivite} implique la conjecture de A.~Escassut et M.C. Sarmant, pour les fonctions holomorphes.
\begin{exemple}\label{exemple}
	Consid\'erons le polyn\^ome $P(z) = z + z^2 \in \CC_p[z]$.
	La boule unit\'e ouverte $B = \{ |z| < 1 \}$ est fix\'ee par $P$ et la restriction de $P$ \`a $B$ pr\'eserve la distance induite par la norme $| \cdot |$.
	En particulier $P|_B$ satisfait la propri\'et\'e {\rm ($\ref{eqn}$)}.

	Clairement $P|_B$ est un \'el\'ement analytique.
	On montrera que $P|_B$ n'est pas bi-analytique.
	On utilisera la notion de {\rm bonne r\'eduction} (introduite dans \cite{MS}) et la notion de {\rm r\'eduction non-triviale}~; voir \cite{hyp} Section~5.1.

	Supposons par contradiction que l'inverse $g : B \longrightarrow B$ de $P|_B$ est un \'el\'ement analytique.
	Il existe alors une fonction rationnelle $R \in \CC_p(z)$ sans p\^oles dans $B$ et telle que
\begin{equation}\label{estimation}
	\sup \{ |R(z) - g(z)| \mid z \in B \} < 1.
\end{equation}

	Comme $g(B) = B$ on a aussi $R(B) = B$ et donc $R$ a une r\'eduction non-triviale $\tilde{R} \in \tilde{\CC}_p(z)$, o\`u $\tilde{\CC}_p$ d\'esigne le corps r\'esiduel de $\CC_p$.
	
	D'autre part, le polyn\^ome $P$ a une bonne r\'eduction, \'egale \`a  $\tilde{P}(z) = z + z^2 \in \tilde{\CC}_p[z]$.
	Donc $\deg(\tilde{R} \circ \tilde{P}) > 1$ et par cons\'equent la fonction rationnelle $R \circ P(z) - z \in \CC_p(z)$ a une r\'eduction non-triviale, \'egale \`a $\tilde{R} \circ \tilde{P}(z) - z \in \tilde{\CC}_p(z)$.
	En particulier
$$
	\sup \{ |R \circ P(z) - z| \mid z \in B \} = 1.
$$
	On obtient une contradiction avec {\rm ($\ref{estimation}$)}.
\end{exemple} 
%%%%%%%%%%%%%%%%%%%%%%%%%%%%%%%%%%%%%%%%%%%%%%%
\section{Une r\'eduction.}
%%%%%%%%%%Preuve du Th\'eor\`eme~\ref{injectivite}.}
%%%%	Notons que la propri\'et\'e ($\ref{eqn}$) reste vraie apr\`es post-composition par une homographie.
	Maintenant on r\'eduit le Th\'eor\`eme~\ref{injectivite} \`a la proposition suivante.

	Notons que toute fonction qui pr\'eserve la distance induite par la norme $| \cdot |$, satisfait la propri\'et\'e ($\ref{eqn}$).
\begin{proposition}\label{proposition}
	Soit $Y \subset \CC_p$ un affino\"{\i}de et soit $f$ une fonction holomorphe d\'efinie sur $Y$.
	Si $f$ est injective, alors il existe une homographie $\varphi$ telle que $\varphi \circ f (Y) \subset \CC_p$ et telle que $\varphi \circ f$ pr\'eserve la distance induite par la norme $| \cdot |$.
\end{proposition}
\begin{proof}[Preuve du Th\'eor\`eme~\ref{injectivite}.]
	Soit $f$ une fonction holomorphe non-constante d\'efinie sur $X$ et satisfaisant la propri\'et\'e ($\ref{eqn}$).
	Notons que si pour un certain $x \in X$ on a $f'(x) = 0$, alors la propri\'et\'e ($\ref{eqn}$) implique que $f$ est constante \'egale \`a $f(x)$.
	On conclut que $f'(x) \neq 0$ pour tout $x \in X$.
	Par cons\'equent si $x, y \in X$ sont distincts, on a
$$
	|f(x) - f(y)|^2 = |x- y|^2 |f'(x) \cdot f'(y)| \neq 0,
$$
et donc $f(x) \neq f(y)$.

	Supposons d'autre part que $f$ est une fonction holomorphe et injective, d\'efinie sur $X$.
	Etant donn\'es $x, y \in X$ distincts, soit $Y \subset X$ un affino\"{\i}de contenant $x$ et $y$.
	Par la Proposition~\ref{proposition}, il existe une homographie $\varphi$ telle que $\varphi \circ f (Y) \subset \CC_p$ et telle que $\varphi \circ f|_Y$ pr\'eserve la distance induite par la norme $| \cdot |$.

	En particulier la fonction $\varphi \circ f |_Y$ satisfait la propri\'et\'e ($\ref{eqn}$).
	Comme $\varphi$ est une homographie, ceci implique que la fonction $f$ satisfait aussi la propri\'et\'e ($\ref{eqn}$) pour $x, y \in X$.
%%%%%	Comme on a choisit $x, y \in X$ de fa\c{c}on arbitraire, on conclut que $f$ satisfait la propri\'et\'e ($\ref{eqn}$).
\end{proof}
%%%%%%%%%%%%%%%%%%%%%%%%%%%%%%%%%%%%%%%%%%%%%%%%%%%
\section{Espace hyperbolique $p$-adique.}
	Pour des r\'ef\'erences \`a cette section on renvoie le lecteur \`a \cite{Y} ou \cite{hyp}.
	
	Consid\'erons la relation d'\'equivalence $\sim$ sur $\CC_p \times \RR$ d\'efinie par
$$
	(w, r) \sim (w', r') \mbox{ si et seulement si $r = r'$ et $p^r \ge |w - w'|$}.
$$
	On consid\`ere la distance $d$ sur $\CC_p \times \RR$ induite par la distance usuelle de $\RR$~: la distance entre les points repr\'esent\'es par $(w, r)$ et $(w', r')$ est donn\'ee par
$|r - r'|$ si $|w - w'| \le p^{\max\{r, r'\} }$ et en g\'en\'eral par
$$
	2 \max \{r, r', \log_p|w - w'| \} - r - r'.
$$
	L'{\it espace hyperbolique} $\HH_p$ est par d\'efinition l'espace m\'etrique $(\CC_p \times \RR/ \sim , d)$ (voir Appendice~2 de \cite{hyp}, o\`u on a consid\'er\'e la compl\'etion de cet espace m\'etrique).
	Les {\it points rationnels} de $\HH_p$ sont par d\'efintion les points repr\'esent\'es par $(w, r)$, avec $r$ rationnel.

	Le groupe des homographies agit par isom\'etries sur $\HH_p$~;
	pour une homographie $\varphi$, on note $\varphi_*$ l'action sur $\HH_p$ correspondante.
	Cette action est d\'etermin\'ee, au niveau de $\CC_p \times \RR$, par les d\'efinitions suivantes~:
\begin{itemize}
	\item Si $T(w) = w + t$, $t \in \CC_p$,  alors $T_*((w, r)) = (w + t, r)$.
	\item Si $M(w) = \lambda w$, $\lambda \in \CC_p^*$, alors $M_*((w, r)) = (\lambda w, r + \log_p \lambda)$.
	\item Si $I(w) = 1/w$, alors $I_*((w, r)) = (1/w, - r)$.
\end{itemize}
	Notons que cette action du groupe des homographies est transitive sur les points rationnels de $\HH_p$.

\

	Etant donn\'e $r \in \RR$, soit $\cs_r$ le point de $\HH_p$ repr\'esent\'e par $(0, r)$.
	On note $(0, \cs_r) = \{ \cs_{r'} \mid r' < r \}$ et $(0, \cs_r] = \{ \cs_{r'} \mid r' \le r \}$, et pour $r_0 < r_1$ on note $[\cs_{r_0}, \cs_{r_1}] = \{ \cs_r \mid r_0 \le r \le r_1 \}$.

	Etant donn\'e $z \in \pp = \CC_p \cup \{ \infty \}$ et $\cs \in \HH_p$ on peut trouver $r \in \RR$ et une homographie $\varphi$, tel que $\varphi(z) = 0$ et $\varphi_*(\cs) = \cs_r$.
	Alors on d\'efinit $(z, \cs) = \varphi_*^{-1}((0, \cs_r))$ et $(z, \cs] = \varphi_*^{-1}((0, \cs_r])$, qui ne d\'ependent pas du choix de $\varphi$.
	On appelle $(z, \cs)$ et $(z, \cs]$ des {\it demi-g\'eod\'esiques}.
	Pour $\cs_0, \cs_1 \in \HH_p$ on d\'efinit $[\cs_0, \cs_1]$ de fa\c{c}on analogue.
%%%%%%%%%%%%%%%%%%%%%%%%%%%%%%%%%%%%%%%%%%%%%%%%%%%%%%%%%
\section{Preuve de la Proposition~\ref{proposition}.}
	Soit $Y = B - (B_1 \cup \ldots \cup B_n)$ un affino\"{\i}de.
	Consid\'erons $a \in B$ et posons $r = \diam(B) \in |\CC_p^*|$.
	On note $\cs$ le point (rationnel) de $\HH_p$ repr\'esent\'e par $(a, \log_p r) \in \CC_p \times \RR$.
	On appelle $\widehat{Y} = \cup_{z \in Y} (z, \cs]$ l'{\it enveloppe convexe de} $Y$~; voir \cite{hyp} Section~3.3.

	Toute fonction holomorphe $f$ d\'efinie sur $Y$ induit une application
$
	f_* : \widehat{Y} \longrightarrow \HH_p
$, voir \cite{Y}.
De plus, si $f$ est injective, alors $f_*$ induit une isom\'etrie entre $\widehat{Y}$ et $f_*(\widehat{Y})$.
	
\

	La Proposition~\ref{proposition} est une cons\'equence imm\'ediate de la proposition suivante.
\begin{proposition}
	On garde les notations pr\'ec\'edentes.
	Alors on a les propri\'et\'es suivantes.
	\begin{enumerate}
		\item
			Il existe une homographie $\varphi$ telle que $(\varphi \circ f)(Y) \subset B$ et $(\varphi \circ f)_*(\cs) = \cs$.
		\item
			Si $f(Y) \subset B$ et $f_*(\cs) = \cs$, alors $f$ pr\'eserve la distance induite par la norme $| \cdot |$.
	\end{enumerate}
\end{proposition}
\begin{proof}[Preuve.]

\noindent
$1.-$	Apr\`es un changement de coordonn\'ee affine au d\'epart on suppose $B = \{ |z| \le 1 \}$.
	Comme le groupe des homographies agit de fa\c{c}on transitive sur les points rationnels de $\HH_p$, on peut trouver une homographie $\varphi$ telle que $(\varphi \circ f)_* (\cs) = \cs$.

	Soit $\co_p = B = \{|z| \le 1 \}$ l'anneau des entiers de $\CC_p$ et soit ${\mathfrak m}_p = \{ |z| < 1 \}$ son id\'eal maximal.
	On \'etend la projection de $\co_p$ \`a $\tilde{\CC}_p = \co_p/{\frak m}_p$, \`a une projection $\pi$ de $\pp$ \`a ${\mathbb P}(\tilde{\CC}_p)$, par $\pi^{-1}(\infty) = \pp - \co_p$.
	Notons que deux points $z, z' \in \pp$ ont la m\^eme projection dans ${\mathbb P}(\tilde{\CC}_p)$ si et seulement si l'intersection $(z, \cs) \cap (z', \cs)$ est non-vide.

	Comme $(\varphi \circ f)_*$ induit une isom\'etrie entre $\widehat{Y}$ et $(\varphi \circ f)_* (\widehat{Y})$, pour tout $z \in Y$ on a $(\varphi \circ f)_*((z, \cs)) = (\varphi \circ f (z), \cs)$.
	Par cons\'equent $\varphi \circ f$ induit une application de $\pi(Y) \subset {\mathbb P}(\tilde{\CC}_p)$ \`a ${\mathbb P}(\tilde{\CC}_p)$.
	On sait que cette application est une restiction d'une fonction rationnelle $\tilde{f} \in \tilde{\CC}_p(z)$, voir \cite{Y}.
	De plus, comme $f$ est injective sur $Y$, le degr\'e de $\tilde{f}$ est \'egale \`a 1.
	Donc on peut choisir $\varphi$ de telle fa\c{c}on que $\tilde{f}$ soit \'egale \`a l'identit\'e.
	Comme $Y \subset B = \co_p$, dans ce cas on a $\varphi \circ f(Y) \subset \co_p = B$.

\

$2.-$	Etant donn\'e deux points $x, y \in B$, soit $\cs_{x, y}$ le point de $\HH_p$ d\'etermin\'e par
$(x, \cs] \cap (y, \cs] = [\cs_{x, y}, \cs]$.
	On a
$$
	|x - y| = \diam(B) \cdot p^{-d(\cs_{x, y}, \cs)}.
$$

	Comme $f_*(\cs) = \cs$ et $f_*$ induit une isom\'etrie entre $\widehat{Y}$ et $f_*(\widehat{Y})$, pour tout $z \in Y$ l'image de la demi-g\'eod\'esique $(z, \cs]$ par $f_*$ est la demi-g\'eod\'esique $(f(z), \cs]$.
	Donc pour $x, y \in Y$ distincts on a $f_*(\cs_{x, y}) = \cs_{f(x), f(y)}$ et par cons\'equent
$$
	d(\cs_{f(x), f(y)}, \cs) = d(f_*(\cs_{x, y}), f_*(\cs)) = d(\cs_{x, y}, \cs).
$$

Comme $f(x), f(y) \in B$ on obtient
$$
	|f(x) - f(y)| = \diam(B) \cdot p^{- d(\cs_{f(x), f(y)}, \cs)}
			= \diam(B) \cdot p^{- d(\cs_{x, y}, \cs)} = |x - y|.
$$
\end{proof}
%%%%%%%%%%%%%%%%%%%%%%%%%%%%%%%%%%%%%%%%%%%%%%%%%%%%%%%%%%%%
%%%  Remerciements  %%%
%%%%%%%%%%%%%%%%%%%%%%%
\section*{Remerciements.}
{\it Je remercie A.~Escassut et J.C.~Yoccoz pour des discussions relies \`a ce travail.
	Mes remerciments vont aussi \`a M.~Flexor qui a fait des corrections d'orthographie.
	Ce papier a \'et\'e ecrit pendent le ``Research Trimester on Dynamical Systems'' \`a Pisa.}

\bibliographystyle{plain}

\begin{thebibliography}{R-L2}

\bibitem{Es95} A. Escassut. {\it Analytic elements in p-adic analysis.}
World Scientific Publishing, $1995$.

\bibitem[ES]{ES} A. Escassut, M.C. Sarmant. {\it Injectivity, Mittag-Leffler series and Motzkin products.}
Ann. Sci. Math. Qu\'ebec, {\bf 16} (1992), 155-173.

\bibitem[FvP]{FvP} J. Fresnel, M. van der Put. {\it G\'eom\'etrie Analytique rigide et applications}. PM 18, Birkh$\ddot{a}$user $1981$.

\bibitem[MS]{MS} P. Morton, J. Silverman. {\it Periodic points, multiplicities, and dynamical units}. J. Reine Agnew. Math. {\bf 461} $(1995), \, 81-122$.

\bibitem[Mo]{Mo} E. Motzkin. {\it L'arbre d'un quasi connexe : un invariant conforme $p$-adique}. Groupe d'\'etude d'analyse ultram\'etrique, $9$eme ann\'ee, {\bf 3} $(1981/82), \, 18$p.

\bibitem[R-L1]{these} J. Rivera-Letelier. {\it Dynamique des fonctions rationnelles sur des corps locaux.} Th\`ese, Orsay $2000$.

\bibitem[R-L2]{hyp} J. Rivera-Letelier. {\it Espace hyperbolique $p$-adique et dynamique des fonctions rationnelles.}
	A para\^{\i}tre dans \,  Comp. Math.

\bibitem[R-L3]{elements} J. Rivera-Letelier. {\it Une caract\'erisation locale des \'el\'ements bi-analytiques.} Texte 2002.

\bibitem[Y]{Y} J.C. Yoccoz. Notes d'un cours au Coll\`ege de France, 2001/2002.

\end{thebibliography}

\end{document}